\documentclass[12pt]{article}
\usepackage{theorem,amsfonts,amssymb}
\newtheorem{theorem}{Theorem}[section]
\newtheorem{proposition}{Proposition}[section]
\newtheorem{lemma}{Lemma}[section]
\newtheorem{corollary}{Corollary}[section]
\theorembodyfont{\upshape}
\newtheorem{definition}{Definition}
\newtheorem{remark}{Remark}[section]
\newtheorem{example}{Example}[section]
\newtheorem{proof}{Proof}

\newtheorem{acknowledgement}{Acknowledgement}

\newcommand{\bt}{\begin{theorem}}
\newcommand{\et}{\end{theorem}}
\newcommand{\bl}{\begin{lemma}}
\newcommand{\el}{\end{lemma}}
\newcommand{\bp}{\begin{proposition}}
\newcommand{\ep}{\end{proposition}}
\newcommand{\bex}{\begin{example}}
\newcommand{\eex}{\end{example}}
\newcommand{\bc}{\begin{corollary}}
\newcommand{\ec}{\end{corollary}}
\newcommand{\bo}{\begin{proof}}
\newcommand{\eo}{\end{proof}}
\newcommand{\bd}{\begin{definition}}
\newcommand{\ed}{\end{definition}}
\newcommand{\br}{\begin{remark}}
\newcommand{\er}{\end{remark}}
\newcommand{\be}{\begin{enumerate}}
\newcommand{\ee}{\end{enumerate}}

\begin{document}

\title{A stochastic difference equation with stationary noise on groups}
\author{C. R. E. Raja}
\date{}
\maketitle

\let\ol=\overline
\let\epsi=\epsilon
\let\vepsi=\varepsilon
\let\lam=\lambda
\let\Lam=\Lambda
\let\ap=\alpha
\let\vp=\varphi
\let\ra=\rightarrow
\let\Ra=\Rightarrow
\let \Llra=\Longleftrightarrow
\let\Lla=\Longleftarrow
\let\lra=\longrightarrow
\let\Lra=\Longrightarrow
\let\ba=\beta
\let\ov=\overline
\let\ga=\gamma
\let\Ba=\Delta
\let\Ga=\Gamma
\let\Da=\Delta
\let\Oa=\Omega
\let\Lam=\Lambda
\let\un=\upsilon

\newcommand{\f}{{\rm frac}}

\newcommand{\Ad}{{\rm Ad}}
\newcommand{\cK}{{\cal K}}
\newcommand{\pc}{T_{si}}
\newcommand{\Aut}{{\rm Aut}}
\newcommand{\cR}{{\cal R}}
\newcommand{\Spr}{{\rm Spr}}
\newcommand{\cG}{{\cal G}}
\newcommand{\cH}{{\cal H}}
\newcommand{\G}{{\mathbb G}}
\newcommand{\Z}{{\mathbb Z}}
\newcommand{\Q}{{\mathbb Q}}
\newcommand{\cN}{{\cal N}}
\newcommand{\cS}{{\cal S}}
\newcommand{\N}{{\mathbb N}}
\newcommand{\R}{{\mathbb R}}
\newcommand{\C}{{\mathbb C}}
\newcommand{\T}{{\mathbb T}}

\begin{abstract}
We consider the stochastic difference equation $$\eta _k = \xi _k
\phi (\eta _{k-1}), ~~~~  k \in \Z $$ on a locally compact group $G$
where $\xi _k$ are given $G$-valued random variables, $\eta _k$ are
unknown $G$-valued random variables and $\phi$ is an automorphism of
$G$. This equation was considered by Tsirelson and Yor on
one-dimensional torus.  We consider the case when $\xi _k$ have a
common law $\mu$ and prove that if $G$ is a pointwise distal group
and $\phi$ is a distal automorphism of $G$ and if the equation has a
solution, then extremal solutions of the equation are in one-one
correspondence with points on the coset space $K\backslash G$ for
some compact subgroup $K$ of $G$ such that $\mu$ is supported on
$Kz= z\phi (K)$ for any $z$ in the support of $\mu$. We also provide
a necessary and sufficient condition for the existence of solutions
to the equation.
\end{abstract}



\begin{section}{Introduction}

Let $G$ be a locally compact group.  Consider the stochastic
difference equation on $G$
$$\eta _k = \xi _k \phi (\eta _{k-1}) , ~~~~ k \in -\N \eqno (Y) $$
where $\eta_k$ and $\xi _k$ are $G$-valued random variables and
$\phi$ is an automorphims of $G$.  The random variables $(\xi _k)$
are given and are called the noise process of the equation $(Y)$. We
are interested in finding the law of the unknown process $(\eta
_k)$. We further assume that for $k \in \Z$, the random variable
$\xi_k$ is independent of $\eta _j$ for $j<k$ and this assumption
will be enforced whenever equation of type $(Y)$ is considered.

B. Tsirelson \cite{Ts} considered the following stochastic
difference equation on the real line
$$\eta _k = \xi _k +\f (\eta _{k-1} ) ~~~~ k \in -\N \eqno (T_1)$$
where $\f (x)$ is the fractional part of $x \in \R$ and $(\xi_k)$ is
a given stationary Gaussian noise process in order to obtain his
celebrated example of the stochastic differential equation
$$dX_t = dB_t+ b^+(t, X)dt, ~~ X(0)=0 \eqno (T_2)$$
which has unique solution that is not strong where $(B_t)$ is the
one-dimensional Brownian motion.  It was also noted that under some
conditions solution of stochastic difference equation $(T_1)$
determines solution of Tsirelson's stochastic differential equation
$(T_2)$ (see \cite{Ts} for more details).

It is easy to see that the set of all solutions $(\eta _k)$ of
equation $(Y)$ is a convex set, hence by extremal solution we mean
the extreme point of the convex set of all solutions.

M. Yor \cite{Yo} formulated equation $(T_1)$ in the form of equation
$(Y)$ on the one-dimensional torus $\R /\Z$ when $\phi$ is the
identity automorphism and $(\xi _k)$ is a general noise process.  In
particular, \cite{Yo} proved that extremal solutions of the equation
$(Y)$ are in one-one correspondence with points on the coset space
$(\R/\Z)/M$ where $M$ is a closed subgroup of $\R/\Z$. When $\phi$
is the identity automorphism, equation $(Y)$ was considered on
general compact groups \cite{AUY} and when the noise law $(\xi _k)$
is stationary, equation $(Y)$ was considered on abelian groups
\cite{Ta}.  The main results of \cite{AUY} amd \cite{Ta} extend the
result of \cite{Yo} and proved that the extremal solutions can be
identified with $G/H$ where $H$ is a certain compact subgroup of
$G$. Assuming that the noise process $(\xi _k)$ is stationary, we
extend the result of \cite{Yo} to a larger class of groups called
pointwise distal groups (that is, $e$ is not a limit point of $\{g^n
xg^{-n} \mid n \in \Z \}$ for any $g,x\in G$ and $x\not =e$) that
includes nilpotent groups, compact groups, discrete groups and
connected groups of polynomial growth and when the automorphism
$\phi$ is distal (that is, $e$ is not the limit of $\{ \phi ^n(x)
\mid n \in \Z \}$ for any $x\in G$ and $x\not =e$).

\bt\label{mt} Let $G$ be a locally compact pointwise distal group
and $\phi$ be a distal automorphism of $G$.  Let $(\xi _k)_{k\in
\Z}$ be $G$-valued random variables with common law $\mu$. Suppose
the equation
$$\eta _k = \xi _k \phi (\eta _{k-1} ), ~~~~k\in \Z \eqno (1) $$
has a solution.  Then there exists a compact subgroup $K_\mu$ such
that for any $z$ in the support of $\mu$, $\mu z^{-1}$ is supported
on $K_\mu = z \phi (K_\mu )z^{-1}$ and we have a one-one
correspondence between left $K_\mu$-invariant probability measures
$\lam$ on $G$ and the laws $(\lam _k)$ of the solutions $(\eta _k)$
of the equation $(1)$ given by
$$\lam _k = \left \{ \begin{array}{ll}
z\phi (z) \cdots \phi ^{k-1}(z) \phi ^k (\lam ) & k\geq 1 \\
\lam & k=0 \\ \phi ^{-1} (z^{-1})\cdots \phi ^{k}(z^{-1}) \phi
^k(\lam )& k<0.
\end{array}\right. $$ Moreover, extremal solutions $(\eta _k)$ of the
equation (1) are in one-one correspondence with the elements of the
coset space $K_\mu\backslash G$.\et

\end{section}

\begin{section}{Preliminaries}

Let $G$ be a locally compact group and $\phi$ be an automorphism of
$G$.  For a (regular Borel) probability measure $\mu$ on $G$, we
define a probability measures $\check \mu$ and $\phi (\mu )$ on $G$
by $\check \mu (E) = \mu (E^{-1})$ and $\phi (\mu )(E) = \mu (\phi
^{-1}(E))$ for any Borel subset $E$ of $G$.

For any two probability measures $\mu$ and $\lam$, the convolution
of $\mu$ and $\lam$ is denoted by $\mu *\lam$ and is defined by
$$\mu *\lam (E) = \int \mu (Ex^{-1}) d\lam (x)$$ for any Borel subset $E$
of $G$.  For $n \geq 1$ and for a probability measure $\mu$ on $G$,
$\mu ^n$ denotes the $n$-th convolution power of $\mu$.

For $x \in G$ and a probability measure $\mu$ on $G$, $x\mu$ (resp.,
$\mu x$) denotes $\delta _x *\mu$ (resp., $\mu *\delta _x$).

For a compact subgroup $K$ of $G$, $\omega _K$ denotes the
normalized Haar measure on $K$ and a probability measure $\lam$ on
$G$ is called left $K$-invariant if $x\lam = \lam $ for all $x \in
K$ (which is equivalent to $\omega _K*\lam = \lam$, by Theorem 1.2.7
of \cite{He}).

We say that a sequence $(\lam _n)$ of probability measures on $G$
converges (in the weak topology) to a probability measure $\lam $ on
$G$ if $\int f d\lam _n \to \int f d\lam$ for all continuous bounded
functions on $G$.

A set $\cal F$ of probability measures on $G$ is said to be
uniformly tight if for each $\epsi >0$ there is a compact set
$C_\epsi$ of $G$ such that $\rho (C_\epsi )>1-\epsi$ for all $\rho
\in {\cal F}$.  It follows from Prohorov's Theorem that ${\cal F}$
is uniformly tight if and only if $\cal F$ is relatively compact in
the space of probability measures on $G$ equipped with weak topology
(cf. Theorem 1.1.11 of \cite{He}).

Let $(\xi _k)_{k\in \Z}$ be $G$-valued random variables. We are
interested in investigating the laws of random variables $(\eta _k)$
that satisfies the stochastic difference equation
$$\eta _k = \xi _k \phi (\eta _{k-1} ), ~~~~ k \in\Z \eqno (1)$$ and
$\xi _k$ is independent of $\eta _j$ for all $j <k$: \cite{AUY} and
\cite{Yo} considered only negative $k$, we also could have
considered $k\in -\N$ but that would not have made any difference.
Since we are interested in the law of the solutions of equation
$(1)$, we would be studying the corresponding convolution equation

$$\lam _k = \mu _k * \phi (\lam _{k-1}) \eqno (1')$$ for all $k \in \Z$
where $\mu _k$ and $\lam _k$ are the laws of $\xi _k$ and $\eta _k$
respectively.

We consider equation (1) when $\xi _k$ is stationary on groups of
the following type:

\bd A locally compact group $G$ is called pointwise distal if for
any $x,g \in G\setminus (e)$, $e$ is not in the closure of $\{ g^n x
g^{-n} \mid n \in \Z \}$. \ed

It is easy to see that discrete groups, in general groups having
small invariant neighborhoods are pointwise distal.  The class of
pointwise distal groups includes all compact extensions of nilpotent
groups, connected groups of polynomial growth (cf. \cite{Ro}),
p-adic Lie groups of type $R$ (cf. \cite{Ra1}).  The class of
pointwise distal groups appear in the study of shifted convolution
powers which is crucial to our study (cf. \cite{RS}).

The following type of automorphisms will be considered.

\bd An automorphism $\phi$ of a locally compact group $G$ is called
distal if $e$ is not in the closure of $\{ \phi ^n (x) \mid n \in \Z
\}$ for any $x \in G\setminus (e)$.  \ed

It is easy to see that if $\phi$ preserves a metric on $G$, then
$\phi$ is distal on $G$.  All unipotent matrices on
finite-dimensional vector spaces are distal.  It can easily be
verified that inner automorphisms of pointwise distal groups are
distal.  In the final section we give examples of compact groups in
which all automorphisms are distal.

Given a automorphism $\phi$ on a locally compact group $G$, the
following type of group is useful in our approach to the equation
(1): the semidirect product of $\Z$ and $G$, denoted by $\Z \ltimes
_\phi G$ is defined by $$(n , g)(m ,h) = (n+m, g\phi ^n(h))$$ for
all $n, m \in \Z$ and $g, h\in G$.  It is easy to see that $\Z
\ltimes _\phi G$ is a locally compact group with $G$ as an open
subgroup.

Given a probability measure $\mu$ on a locally compact group $G$ and
an automorphisms $\phi$ of $G$, we will also be studying the
probability measure $1\otimes \mu$ on $\Z \ltimes _\phi G$ defined
by $$1\otimes \mu (A\times B) = \delta _1(A) \mu (B)$$ for any
subset $A$ of $\Z$ and any Borel subset $B$ of $G$.
\end{section}

\begin{section}{Dissipating measures}

A probability measure $\mu$ on a locally compact group $G$ is called
dissipating if for any compact set $C$ in $G$, $\sup _{x\in G}\mu ^n
(xC)\to 0$.

In the study of dissipating measures the smallest closed normal
subgroup a coset of which contains the support of $\mu$ plays a
crucial role, let $N_\mu$ denote this normal subgroup of $G$. Let
$G_\mu$ be the closed subgroup generated by the support of $\mu$. If
$G_\mu$ is non-compact and $G_\mu/N_\mu$ is compact, then \cite{JRW}
showed that $\mu$ is dissipating.  \cite{Ja1} and \cite{Ja2}
provided many sufficient conditions (on $G_\mu$ or on $\mu$) for
$\mu$ to be dissipating, for instance if $G_\mu$ is not amenable,
then $\mu$ is dissipating (Corollary 1.6 of \cite{Ja1} or Corollary
3.6 of \cite{Ja2}).

We will now provide a necessary and sufficient condition so that the
equation $(1)$ has a solution.

\bp\label{pr1}Let $G$ be a locally compact group and $(\xi _k)_{k\in
\Z}$ be $G$-valued random variables with common law $\mu$. Let
$\phi$ be an automorphism of $G$.  Then there is a solution $(\eta
_k)$ of the equation
$$ \eta _k = \xi _k \phi (\eta _{k-1}), ~~~~ k \in \Z$$ if and only if
the probability measure $\rho = 1\otimes \mu$ on $\tilde G = \Z
\ltimes _\phi G$ is not dissipating. \ep

\br If $\phi$ is trivial in Proposition \ref{pr1}, then $\rho ^n =
n\otimes \mu^n$, hence for any compact set $C$ of $G$, $\sup _{a\in
\tilde G} \rho ^n(Ca) = \sup _{x\in G}\mu ^n(Cx)$. Since $G$ is open
in $\tilde G$, $\rho$ is dissipating if and only if $\mu$ is
dissipating. Thus, when $\phi$ is trivial, the equation in
Proposition \ref{pr1} has a solution if and only if $\mu$ is
dissipating. \er

\bo If there is a sequence $(\eta _k)$ of $G$-valued random
variables such that
$$\eta _k = \xi _k \phi (\eta _{k-1} )$$ for all $k \in \Z$. Let
$\lam _k$ be the law of $\eta _k$.  Then the corresponding
convolution equation is
$$\lam _k = \mu
* \phi (\lam _{k-1}) $$ for all $k \in \Z$.  Iterating the
convolution equation we get that
$$\lam _k = \mu *\phi (\mu ) *\cdots *\phi ^{j-1}(\mu )
* \phi ^j(\lam _{k-j} )$$ for all $k \in Z$ and all $j \geq 1$. It follows
from Theorems 1.2.21 (iii) of \cite{He} that there is a sequence
$(g_j)$ in $G$ such that the sequence $(\mu *\phi (\mu )
*\cdots *\phi ^{j-1}(\mu ) g_j^{-1})$ is uniformly tight.

Consider the group $\tilde G = \Z \ltimes _\phi G$ and the measure
$\rho = 1\otimes \mu$ on $\tilde G$.  Then $\rho ^j = j\otimes \mu
*\phi (\mu ) *\cdots *\phi ^{j-1}(\mu )$, hence we get that $\rho ^j
(j, g_j)^{-1} = j\otimes \mu *\phi (\mu ) *\cdots *\phi ^{j-1}(\mu )
(-j, \phi ^{-j}(g_j^{-1}))= \mu *\phi (\mu )
*\cdots *\phi ^{j-1}(\mu ) g_j^{-1}$.  This implies that $(\rho
^j(j,g_j)^{-1})$ is uniformly tight. Then there is a compact set $C$
of $\tilde G$ such that $\rho ^j (C(j,g_j)) > {1\over 2} $ for all
$j \geq 1$. This proves that $\rho$ is not dissipating.

For the converse, suppose $\rho$ on $\tilde G = \Z \ltimes _\phi G$
is not dissipating. We first assume that $G$ is separable.  Then by
Theorem 3.1 of \cite{Cs}, there is a sequence $(n_j, g_j)$ in
$\tilde G$ such that $(\rho ^j (n_j, g_j)^{-1})$ converges (see also
\cite{To}). Now since $\rho ^j = j\otimes \mu
*\phi (\mu )*\cdots * \phi ^{j-1}(\mu )$, we get that $\rho ^j(n_j
,g_j)^{-1} = (j-n_j)\otimes \mu *\phi (\mu )*\cdots
* \phi ^{j-1}(\mu )\phi ^{j-n_j}(g_j^{-1})$, hence $(\mu *\phi (\mu
)*\cdots * \phi ^{j-1}(\mu )x_j)$ converges for $x_j = \phi
^{j-n_j}(g_j^{-1})$.  Let $\ga = \lim \mu *\phi (\mu )*\cdots * \phi
^{j-1}(\mu )x_j$.  Then $\ga = \mu *\phi (\lim (\mu *\phi (\mu
)*\cdots * \phi ^{j-2}(\mu )x_{j-1})* \phi(x_{j-1}^{-1})x_j$, hence
by Theorems 1.2.21 and 1.1.11 of \cite{He}, we get that $(\phi
(x_{j-1}^{-1})x_j)$ is relatively compact.  If $x$ is a limit point
of $(\phi (x_{j-1}^{-1})x_j)$, then $\ga = \mu *\phi (\ga ) x$.

Let $\lam _0 =\ga$ and for $k \geq 1$, let $\lam _k = \ga x^{-1}\phi
(x^{-1})\cdots \phi ^{k-1}(x^{-1})$.  Then for $k\geq 0$, $\lam
_{k+1} = \lam _k \phi ^k(x^{-1})$.  If $\lam _k  = \mu *\phi (\lam
_{k-1})$ for $k\geq 1$, then $\lam _{k+1} = \mu *\phi (\lam
_{k-1})\phi ^{k}(x^{-1})= \mu *\phi (\lam _{k-1}\phi ^{k-1}(x^{-1})
)=\mu *\phi (\lam _k)$.  Since $\lam _1 = \ga x^{-1}= \mu *\phi
(\lam _0)$, it follows from induction that $\lam _{k}= \mu *\phi
(\lam _{k-1})$ for all $k \geq 1$.

For $k <0$, let $\lam _k = \ga \phi ^{-1}(x)\cdots \phi ^{k}(x)$.
Then for $k <0$, $\lam _k= \lam _{k+1} \phi ^k(x)$. If $\lam _k =
\mu *\phi (\lam _{k-1})$ for $k\leq 0$, then $\lam _{k-1} = \lam _k
\phi ^{k-1}(x) = \mu *\phi (\lam _{k-1}\phi ^{k-2}(x)) = \mu *\phi
(\lam _{k-2})$. Since $\lam _0 = \ga = \mu *\phi (\ga ) x = \mu
*\phi (\ga \phi ^{-1}(x)) = \mu *\phi (\lam _{-1})$, it follows from
induction that $\lam _k = \mu *\phi (\lam _{k-1})$ for all $k \leq
0$.  Thus, proving $\lam _k = \mu*\phi (\lam _{k-1})$ for all $k \in
\Z$.

In general, if $G$ is any locally compact group, then replacing $G$
by the smallest $\phi$-invariant closed subgroup of $G$ containing
the support of $\mu$ we may assume that $G$ is $\sigma$-compact.
This implies that $\Z \ltimes _\phi G$ is $\sigma$-compact.  Then by
Theorem 8.7 of \cite{HR}, there exists a $\phi$-invariant compact
normal subgroup $K$ of $G$ such that $G/K$ is separable.  Let $\mu'$
be the image of $\mu$ on $G/K$. Since $K$ is compact, $1\otimes
\mu'$ is not dissipating on $\Z \ltimes _\phi G/K$.  By the previous
case, there exists probability measures $\lam _k '$ on $G/K$ such
that $\lam _k' = \mu'
*\phi (\lam _{k-1}')$ for all $k \in \Z$.  It follows from 1.2.15
(iii) of \cite{He} that there exists probability measures $\lam _k$
on $G$ such that $\lam _k*\omega _K = \lam _k$ and image of $\lam
_k$ on $G/K$ is $\lam _k'$.  Since $K$ is $\phi $-invariant, $\mu
*\phi (\lam _{k-1})*\omega _K = \mu *\phi (\lam _{k-1})$ for all $k
\in \Z$. Since both $\lam _k$ and $\mu *\phi (\lam _{k-1})$ are
projected onto the same probability measure $\lam _k '$ on $G/K$, by
Theorem 1.2.15 (iii) of \cite{He} we get that $\lam _k = \mu*\phi(
\lam _{k-1})$ for all $k \in \Z$. \eo

\end{section}

\begin{section}{Shifted convolution property}

A probability measure $\mu$ on a locally compact group $G$ is said
to have shifted convolution property if either $\mu$ is dissipating
or there is a compact subgroup $K$ of $G$ and a $g\in G$ such that
$\mu ^n g^{-n}\to \omega _K$ and $gKg^{-1}=K$. Shifted convolution
property was studied in details in \cite{RS}, it is shown that all
probability measures on a locally compact group $G$ have shifted
convolution property if and only if the group $G$ is pointwise
distal (see Theorem 6.1 of \cite{RS}). We first prove the following
result which provides a sufficient condition for the existence of
large collection of solutions.

\bp\label{nse} Let $G$ be a locally compact group and $\mu$ be a
probability measure on $G$.  If there
is a compact subgroup $K$ of $G$ such that for any $z$ in the
support of $\mu$, $\mu z^{-1}$ is supported on $K= z\phi (K)z^{-1}$,
then for any left $K$-invariant probability measure $\lam$, $(\lam
_k)$ defined by
$$\lam _k = \left \{ \begin{array}{ll}
z\phi (z) \cdots \phi ^{k-1}(z) \phi ^k (\lam ) & k\geq 1 \\
\lam & k=0 \\ \phi ^{-1} (z^{-1})\cdots \phi ^{k}(z^{-1}) \phi
^k(\lam )& k<0
\end{array}\right. $$ is a solution to equation

$$\lam _k=\mu * \phi (\lam _{k-1}), ~~~~ k \in \Z \eqno (1')$$
for any $z$ in the support of $\mu$. \ep

\bo Assume that there is a compact subgroup $K$ of $G$ such that for
any $x$ in the support of $\mu$, $\mu x^{-1}$ is supported on $K=
x\phi (K)x^{-1}$.  Suppose $z$ is in the support of $\mu$.  Let
$\lam$ be a left $K$-invariant probability measure on $G$ and define
$\lam _0 = \lam$. For $k \geq 1$, let $\lam _k =z\phi (z) \cdots
\phi ^{k-1}(z) \phi ^k (\lam )$. Then $\lam _k = z \phi (\lam
_{k-1})$ for all $k \geq 1$.

We first claim that $\lam _k$ is left $K$-invariant for all $k \geq
0$.  Our claim is based on induction.  For $k \geq 1$, if $\lam
_{k-1}$ is left $K$-invariant, then $\phi (\lam _{k-1})$ is left
$\phi (K)$-invariant, hence for $x\in K$, $x\lam _k =  xz \phi (\lam
_{k-1}) = z \phi (\lam _{k-1}) = \lam _k$ as $\phi (K) = z^{-1}Kz$
implies $z^{-1}xz \in \phi (K)$. This proves that $\lam _k$ is left
$K$-invariant, if $\lam _{k-1}$ is left $K$-invariant.  Since $\lam
_0 = \lam $ is left $K$-invariant, induction implies that $\lam _k$
is left $K$-invariant for all $k \geq 0$.

Since $\mu z^{-1}$ is supported on $K = z\phi (K)z^{-1}$, we get
that $z^{-1}\mu $ is supported on $\phi (K)$.  Since $\lam _{k-1}$
is left $K$-invariant, $\mu *\phi (\lam _{k-1})= z \phi (\lam
_{k-1}) = \lam _k$ for all $k \geq 1$.

For $k < 0$, let $\lam _k =\phi ^{-1}(z^{-1})\cdots \phi
^{k}(z^{-1}) \phi ^k(\lam )$.  Then $\lam _k = \phi
^{-1}(z^{-1})\phi ^{-1}(\lam _{k+1})$ for all $k<0$.

We now claim by induction that $\lam _k$ is left $K$-invariant for
all $k\leq 0$. If for $k<0$, $\lam _{k+1}$ is left $K$-invariant,
then since $z^{-1}Kz= \phi (K)$, we have $\phi ^{-1}(K) = \phi
^{-1}(z)K \phi ^{-1} (z^{-1})$, hence for $x \in K$, $x\lam _k = x
\phi ^{-1} (z^{-1})\phi ^{-1}(\lam _{k+1}) = \phi ^{-1}(z^{-1})\phi
^{-1}(\lam _{k+1}) = \lam _k$.  This proves that $\lam _k$ is left
$K$-invariant for all $k \leq 0$.

For $k< 0$, we have $\lam _k = \phi ^{-1}(z^{-1})\phi ^{-1}(\lam
_{k+1})$, hence $\phi (\lam _k) = z^{-1}\lam _{k+1}$. This implies
that for $k \leq 0$, $\mu *\phi ( \lam _{k-1}) = \mu * z^{-1}\lam _k
=\lam _k$ as $\mu z^{-1}$ is supported on $K$.  Thus, $(\lam _k)$ is
a solution to equation $(1')$. \eo

\bp\label{ns} Let $G$ be a locally compact group and $\phi$ be a
automorphism of $G$.  Let $(\xi _k)$ be a sequence of $G$-valued
random variables with common law $\mu$. Suppose the measure
$1\otimes \mu$ has shifted convolution property on $\Z \ltimes _\phi
G$.  Then for any solution $(\eta _k)$ of the equation
$$\eta _k = \xi _k \phi (\eta _{k-1} ), ~~~~k\in \Z \eqno (2) $$ and
for any $z$ in the support of $\mu$, we have

\be

\item [(1)] there is a compact subgroup $K_\mu$ such that $\mu$ is
supported on $K_\mu z= z\phi (K_\mu ) $,

\item [(2)] $\eta _k = z \phi ( \eta _{k-1} )$ in law,

and when equation $(2)$ has a solution

\item [(3)] there is a one-one correspondence between left $K_\mu$-invariant
probability measures $\lam$ on $G$ and the laws $(\lam _k)$ of the
solutions $(\eta _k)$ of the equation $(2)$ given by
$$\lam _k = \left \{ \begin{array}{ll}
z\phi (z) \cdots \phi ^{k-1}(z) \phi ^k (\lam ) & k\geq 1 \\
\lam & k=0 \\ \phi ^{-1} (z^{-1})\cdots \phi ^{k}(z^{-1}) \phi
^k(\lam )& k<0 .
\end{array}\right. $$
\ee \ep

\bo Let $\tilde G = \Z \ltimes _\phi G$ be the semidirect product of
$\Z$ and $G$ where the $\Z$-action is given $\phi$ and $\rho =
\delta _1 \otimes \mu$.  Suppose $(\eta _k)$ is a solution to
equation $(2)$.  Then by Proposition \ref{pr1}, $\rho$ is not
dissipating. Since $\rho$ has shifted convolution property, there is
a compact subgroup $K$ of $\tilde G$ and $g\in \tilde G$ such that
$\rho ^i g^{-i} \to \omega _K$ and $gKg^{-1} =K$.  By Theorem 4.3 of
\cite{Ei} we get that for any $a$ in the support of $\rho$, $\rho ^i
a^{-i}\to \omega _K$ and $aKa^{-1} = K$ (cf. Remark 1.2 \cite{RS}).

Since $\tilde G/ G \simeq \Z$, we get that $K\subset G$.  Let $z$ be
in the support of $\mu$ and $ a= (1, z)$.  Then $a$ is in the
support of $\rho$.  Let $z_i = \phi ^{-1}(z^{-1}) \cdots \phi
^{-i+1} (z^{-1} )\phi ^{-i}(z^{-1} )$ for $i >1$.  Then we get that
$\rho ^i a^{-i}=\mu * \phi (\mu ) * \cdots *\phi ^{i-1}(\mu ) \phi
^i(z_i) \to \omega _K$.

Let $\lam _k$ be the law of $\eta _k$, $k \in \Z$.  We now claim
that for $k \in \Z$, $\lam _k$ is left $K$-invariant. For $k\in \Z$,
$$\lam _k =\mu *\phi (\lam _{k-1}) = \mu
* \phi (\mu ) * \cdots *\phi ^{i-1}(\mu ) \phi ^i(z_i) * \phi ^i(z_i^{-1})
\phi ^i(\lam _{k-i}), ~~ i \geq 1$$ and hence by Theorems 1.2.21
(ii) and 1.1.11 of \cite{He}, $(\phi ^i(z_i^{-1}) \phi ^i (\lam
_{k-i}) )_{i\geq 1}$ is relatively compact. Thus, for any limit
point $\nu$ of $(\phi ^i(z_i^{-1} )\phi ^i(\lam _{k-i}) )$, we get
that $\lam _k = \omega _K \nu$.  Thus, $\lam _k$ is left
$K$-invariant.

For any $a$ in the support of $\rho$, we have $\rho ^i a^{-i} \to
\omega _K$, hence $$\rho \omega _K a^{-1} = \lim \rho \rho
^{i-1}a^{-i+1}a^{-1} = \lim \rho ^i a^{-i} = \omega _K .$$  This
shows that $\rho$ is supported on $aK = Ka$. Now for $a = (1, z)$,
$a^{-1}\rho = \phi ^{-1}(z^{-1})\phi ^{-1}(\mu )$ is supported on
$K$, hence we get that $z^{-1}\mu$ is supported on $\phi (K)$. This
implies that for $k \in \Z$,
$$\lam _k = \mu * \phi (\lam _{k-1}) = z z^{-1}\mu *\phi (\lam _{k-1}) =
z \phi (\lam _{k-1})$$ as $\lam _{k-1}$ is left $K$-invariant. This
proves (2).

Let $\lam = \lam _0$.  Then $\lam $ is left $K$-invariant and for $k
\geq 1$, $\lam _k = z \phi (\lam _{k-1}) = z\phi (z) \cdots \phi
^{k-1}(z) \phi ^k (\lam )$.  For $k <0$, $\lam = z\phi (\lam _{-1})
= \cdots = z\phi (z) \cdots \phi ^{-k-1}(z) \phi ^{-k}(\lam _k)$,
hence $\lam _k = \phi ^{-1}(z^{-1})\cdots \phi ^{k}(z^{-1}) \phi
^{k}(\lam )$.

For any $z$ in the support of $\mu$, $a=(1, z)$ is in the support of
$\rho$, hence $aKa^{-1}=K$.  This implies that $z\phi (K)z^{-1} =
K$, hence $\phi (K) = z^{-1}Kz$.  Since $z^{-1}\mu$ is supported on
$\phi (K)$, we get that $\mu z^{-1} = z (z^{-1}\mu ) z^{-1}$ is
supported on $z\phi (K)z^{-1}=K$.  This shows that the conditions of
Proposition \ref{nse} are satisfied.  Thus, for a left $K$-invariant
measure $\lam $ if we define $(\lam _k)$ as in the proposition, we
get that $\lam _k$ satisfies $\lam _k = \mu *\phi (\lam _{k-1})$ for
all $k \in \Z$.

\eo

\bc \label{ex} Let $G$, $\phi$ and $(\xi _k)$, $\mu$ be as in Proposition
\ref{ns}.  Suppose $1\otimes \mu$ has shifted convolution property
and $1\otimes \mu$ is not dissipating. Then there exists a compact
subgroup $K$ such that extremal solutions $(\eta _k)$ of the
equation
$$\eta _k = \xi _k \phi (\eta _{k-1} ), ~~~~k\in \Z \eqno (2) $$
are in one-one correspondence with the elements of the coset space
$K\backslash G$.\ec

\bo Let $K$ be as in Proposition \ref{ns}.  Then it follows from
Proposition \ref{ns} that left $K$-invariant measures and laws of
solutions to the equation (2) are in one-one correspondence.

If $\lam$ is a left $K$-invariant measure and $z$ is in the support
of $\lam$.  Let $a \in K$ and $U$ be a neighborhood of $az$. Then
$\lam (U)= \lam (a^{-1}U) >0$.  This implies that $az$ is in the
support of $\lam$.  Thus, support of $\lam$ is a union of cosets of
$K$.  If the support of $\lam$ contains more than one coset, then
the corresponding solution to the equation (2) is not extremal. This
proves the corollary. \eo

We have the following converse to Proposition \ref{ns}.

\bp\label{cns} Let $G$ be a locally compact and $\phi$ be an
automorphism.  Let $(\xi _k)$ be $G$-valued random variables with
common law $\mu$.  Suppose the laws of the solutions $(\eta _k)$ to
the equation
$$\eta _k = \xi _k \phi (\eta _{k-1}),~~~~k \in \Z$$ are left
$K$-invariant for some compact subgroup $K$ of $G$ such that $\mu
z^{-1}$ is supported on $K=z\phi (K)z^{-1}$ for any $z$ in the
support of $\mu$.  Then $1\otimes \mu $ on $\Z \ltimes _\phi G$ has
shifted convolution property.
\ep

\bo Let $\rho =1\otimes \mu$.  We first assume that $G$ is
separable. Since the equation has solutions, by Proposition
\ref{pr1}, $1\otimes \mu$ on $\Z\ltimes _\phi G$ is dissipating.  As
in Proposition \ref{pr1}, there is a $x_j \in G$ and a probability
measure $\ga$ on $G$ such that $\mu
*\phi (\mu )*\cdots *\phi ^j(\mu )x_j \to \ga$ and a solution $(\lam
_k)$ with $\lam _0 = \ga$.  Assumption implies that $\ga$ is left
$K$-invariant.

We now claim by induction that $\mu *\phi (\mu )*\cdots *\phi ^j(\mu
)$ is supported on $Kz\phi (z)\cdots \phi ^j(z)$ for all $j \geq 0$.
If $\mu *\phi (\mu )*\cdots *\phi ^j(\mu )$ is supported on $Kz\phi
(z)\cdots \phi ^j(z)$ for some $j \geq 0$, then $\mu *\phi (\mu
)*\cdots *\phi ^j(\mu )*\phi ^{j+1}(\mu )$ is supported on $Kz\phi
(z)\cdots \phi ^j(z)\phi ^{j+1}(K)\phi ^{j+1}(z)$.  Since $Kz= z\phi
(K)$, we get that $\phi ^k(z)\phi ^{k+1}(K) = \phi ^k(K)\phi ^k(z)$
for all $k \geq 0$.  This shows that $\mu *\phi (\mu )*\cdots *\phi
^j(\mu )*\phi ^{j+1}(\mu )$ is supported on $Kz\phi (z)\cdots \phi
^j(z)\phi ^{j+1}(z)$.  Since $\mu$ is supported on $Kz$, claim
follows from induction.

For $j \geq 1$, let $\sigma _j = \mu *\phi (\mu )*\cdots *\phi
^j(\mu )x_j$.  Then $\sigma _j *\check \sigma _j \to \ga *\check
\ga$ and $\sigma _j x_j^{-1}$ is supported on $Kz\phi (z)\cdots \phi
^j(z)$. This implies that $\sigma _j *\check \sigma _j$ is supported
on $K$, hence $\ga *\check \ga$ is supported on $K$. Since $\ga$ is
left $K$-invariant, $\ga *\check \ga$ is also left $K$-invariant and
hence $\ga *\check \ga = \omega _K$.  Now $\rho ^j = j\otimes \sigma
_{j-1}x_{j-1}^{-1}$, hence $\rho ^j *\check \rho ^j = \sigma
_{j-1}*\check \sigma _{j-1}$ for all $j >1$.  This implies that
$\rho ^j *\check \rho ^j \to \omega _K$.  By Theorem 4.3 of
\cite{Ei}, for any $g$ in the support of $\rho$, $\rho ^jg^{-j}\to
\omega _K$. For any $g $ in the support of $\rho = 1\otimes \mu$,
there is a $z$ in the support of $\mu$ such that $g=(1, z)$ and
hence $gKg^{-1} = z\phi (K)z^{-1} =K$.  This proves that $\rho$ has
shifted convolution property.

If $G$ is any locally compact group, replacing $G$ by the smallest
$\phi$-invariant closed subgroup of $G$ containing the support of
$\mu$, we may assume that $G$ is $\sigma$-compact and hence $\Z
\ltimes _\phi G$ is also $\sigma$-compact. Then by Theorem 8.7 of
\cite{HR}, each neighborhood $U$ of $e$ contains a $\phi $-invariant
compact normal subgroup $K_U$ such that $G/K_U$ is separable. It can
easily be verified that the assumptions in the proposition are valid
for $G/K_U$ with $KK_U/K_U$ and the image of $\mu$ on $G/K_U$.  It
follows from the previous case that image of $\mu$ on $G/K_U$ has
shifted convolution property.  By Proposition 2.3 of \cite{RS} we
get that $\mu$ itself has shifted convolution property. \eo

We now extend the results of \cite{AUY}, \cite{Ta} and \cite{Yo} to
all pointwise distal groups when $\xi _k$ is stationary and $\phi$
is distal on $G$.

\bo $\!\!\!\!\!$ {\bf of Theorem \ref{mt}} If $G$ is pointwise
distal and $\phi$ is distal, then the group $\Z \ltimes _\phi G$ is
a pointwise distal group.  By Theorem 6.1 of \cite{RS}, we get that
$1\otimes \mu$ has shifted convolution property.  Now the result
follows from Proposition \ref{ns} and Corollary \ref{ex}. \eo

\br\label{r1} We would like to remark that if $1\otimes \mu$ does
not have shifted convolution property, then the conclusion on
extreme points of the solutions in Theorem \ref{mt} may not be true
even on compact groups. Consider the two dimensional torus $K =
(\R/\Z )^2$ and $\phi$ be an automorphism of $K$ such that $C(\phi
)= \{ x\in K \mid \phi ^n (x) \to e {~~as}~~ n \to \infty \}$ is
dense in $K$, for instance if we take $\phi$ to be $\phi (x, y) =
(x+y+\Z, x+2y+\Z)$ for all $x,y \in \R$, then $C(\phi ) = \{ (t+\Z ,
({1-\sqrt {5}\over 2})t+\Z) \mid t\in \R \} \simeq \R$ is a vector
(nonclosed) subgroup of $K$ and is dense in $K$. Take $\mu$ to be a
probability measure on $K$ such that support of $\mu$ is a compact
subset contained in $C(\phi )$. Since $\phi $ on $C(\phi )$ is
multiplication by ${3-\sqrt{5}\over 2}$, \cite{Za} implies that
there is a probability measure $\rho$ on $C(\phi )$ such that $\mu
* \phi (\mu ) * \cdots *\phi ^i (\mu ) \to \rho$ in the space of
probability measures on $C(\phi )$. This implies that $\rho = \lim
\mu *\phi (\mu ) *\cdots *\phi ^i(\mu ) = \mu *\phi (\lim (\mu
*\cdots *\phi ^{i-1}(\mu ))= \mu *\phi (\rho )$.
Taking $\lam _k = \rho$ for all $k \in \Z$, we get a stationary
solution to the equation $(1)$. Further, assume that $\mu \not =
\delta _x$ for any $x \in K$.  Then $\lam _k=\rho $ are also not
dirac measures. If $x \lam _k =\lam _k$ for some $x \in K$, then
since $\lam _k(C(\phi )) =1$, we get that $\lam _k(C(\phi )) =1 =
x\lam _k(C(\phi ))$, hence $x \in C(\phi )$ as $C(\phi )$ is a
group.  By Theorem 1.2.4 of \cite{He}, $\{ g \in C(\phi ) \mid g\rho
= \rho \}$ is a compact subgroup of $C(\phi )$.  Since $C(\phi )$ is
a vector group, $C(\phi )$ has no non-trivial compact subgroup and
hence $x = e$. Thus, $\lam _k$ is not left invariant for any
nontrivial compact subgroup of $K$. Also in this case $K_\mu = K$.
Hence, the conclusion on the extreme points of solutions in Theorem
\ref{mt} does not hold. \er

\end{section}

\begin{section}{Examples}

We first provide examples of groups for which the group of
automorphisms is compact.

\vskip0.1in

\noindent (i) {\bf Compact p-adic Lie groups:} Let $K$ be a compact
$p$-adic Lie group. Then $\Aut (K)$ is a compact group (see
Corollary 8.35 of \cite{DS} or \cite{Ra0}). The following are
examples of compact $p$-adic Lie groups.

\be

\item [(a)] If $\Q _p$ is the field of $p$-adic numbers with
valuation $|\cdot |_p$, then $\Z _{p^n} = \{ x \in \Q _p \mid |x|_p
\leq p ^{n-1} \}$ is a compact $p$-adic Lie group.

\item [(b)] The group $GL_k(\Z _p)$ of all invertible $k\times
k$-matrices over $\Z_p$.

\item [(c)] Pro-p group of finite rank, that is a totally
disconnected group of finite rank in which every open normal
subgroup has index equal to some power of $p$.

\ee

\vskip0.1in

\noindent (ii) {\bf A Compact abelian group:}  For $n \geq 1$, let
$A_n$ be the group of all $n$-th roots of unity and $A=\cup A_n$.
Then $A$ is a countable abelian group.  Equip $A$ with discrete
topology. Let $K$ be the dual of $A$.  Then $K$ is a compact
(totally disconnected) metrizable group with dual $A$ (see 24.15 of
\cite{HR}).

Let $K_n$ be the closed subgroup of $K$ such that $K/K_n$ is the
dual of $A_n$.  Since $A_n$ is finite, $K/K_n$ is finite. Then $K_n$
is an open subgroup of $K$.  Now, if $x \in \cap K_n$, then $x\in
K_n$ for all $n \geq 1$.  This implies that $a(x)=1$ for all $a \in
A_n$ and for all $n \geq 1$.  Since $A=\cup A_n$, $a(x) =1$ for all
$a \in A$.  Since $A$ is the dual of $K$, $x=e$.  Thus, $\cap K_n
=e$.

Let $\ap$ be an automorphism of $K$ and $\hat \ap$ be the
corresponding dual automorphism on $A$.  Then it is easy to see that
$\hat \ap (A_n) =A_n $ for all $n \geq 1$.  This implies that $\ap
(K_n) =K_n$. This proves that $(K_n)$ is a sequence of arbitrarily
small characteristic open subgroups, hence the group of
automorphisms of $K$ is compact.

\vskip0.2in

\noindent (iii) {\bf All automorphisms are distal but group of
automorphsims is not compact:} Let $\R/\Z$ be the one-dimensional
torus and $K$ be the compact group in (i) or in (ii).  Let $G= \R/\Z
\times K$ be the direct product of $\R/\Z$ and $K$.  Then $G$ is a
compact group.

Let $\tau$ be an automorphism of $G$.  We now claim that there is an
automorphism $\ap$ of $K$ and a character $\chi$ of $K$ such that
$\tau (z, x) = (z^{\pm 1}\chi(x), \ap (x)) $ for all $(z, x)\in G$.
Since the connected component of identity in $G$ is $\R/\Z$, we get
that $\tau (\R/\Z) = \R/Z$, hence $\tau (z, e) = (z^{\pm 1}, e)$ for
all $z\in \R/\Z$. Let $\ap \colon K \to K$ be defined by $\ap (x) =
p(\tau (1,x))$ where $p\colon G \to K$ is the canonical projection
of $G$ onto $K$.  It is easy see that $\ap$ is a continuous
homomorphism. If $\ap (x) =e$, then $p(\tau (1,x) )=e$ and hence
$\tau (1,x) =(z,e)$ for some $z \in \R/\Z$.  But $\tau (z', e) =(z,
e)$ for $z'=z$ or $z'=z^{-1}$. Since $\tau$ is an automorphism,
$(z', e) = (1, x)$, hence $x=e$. This shows that $\ap$ is injective.
For $x \in K$, let $y \in K$ and $z\in \R/\Z$ be such that $\tau (z,
y) = (1, x)$ as $\tau$ is onto. This implies that $\ap (y) = p(\tau
(1,y))=p(\tau (z,y))=x$. This proves that $\ap$ is bijective.
Continuity of $\ap^{-1}$ follows from open mapping theorem as $K$ is
a compact metrizable group (cf. 5.29 of \cite{HR}).  Thus, $\ap$ is
an automorphism of $K$. Let $\chi \colon K \to \R/\Z$ be defined by
$\chi (x) = q(\tau (1,x))$ where $q\colon G \to \R/\Z$ is the
canonical projection of $G$ onto $\R/\Z$. Then $\chi$ is a
continuous homomorphism. For $z\in \R/\Z$ and $x\in K$, $\tau (z,x)
= (z^{\pm 1},e) \tau (1,x) = (z^{\pm 1},e)(\chi (x), \ap (x)) =
(z^{\pm 1}\chi(x), \ap(x))$. This proves the claim.

We now claim that $\tau$ is distal.  Suppose $(1,e)$ is in the
closure of $\{ \tau ^n(z,x) \mid n \in \Z \}$.  Then $e$ is in the
closure of $\{ \ap ^n(x) \mid n \in \Z \}$.  Since the group of
automorphisms of $K$ is compact, $x=e$.  This implies that $\tau ^n
(z, e)=(z^{\pm 1}, e)$ and hence $z=1$.  Thus, $\tau$ is distal. In
fact, one can show that each $\tau$ generates a compact subgroup.

If $K$ is not a finite group, then the group of automorphisms of $G$
is not a compact group as it is homeomorphic to $\{ {\pm 1} \}
\times \hat K \times {\rm Aut}(K)$ where ${\rm Aut}(K)$ is the group
of automorphisms of $K$.

\end{section}

\begin{acknowledgement}
I would like to thank Prof. Riddhi Shah of JNU for suggesting the if
part of Proposition \ref{pr1} and for many useful discussions, in
particular regarding Remark \ref{r1}.
\end{acknowledgement}

\begin{tabular}{l}
C.\ R.\ E.\ Raja \\
Stat-Math Unit \\
Indian Statistical Institute\\
8th Mile Mysore Road \\
Bangalore 560 059, India.\\
e-mail: creraja@isibang.ac.in
\end{tabular}

\end{document}